\documentclass[12pt]{article}
\usepackage{amsmath,amssymb}
\usepackage{graphicx}
\newcommand{\eqdef}{\stackrel{\text{def}}{=}}

\newcommand{\n}{\nonumber \\}

\newcommand{\ignore}[1]{}
\numberwithin{equation}{section}
\newcommand{\Romannumeral}[1]{\uppercase\expandafter{\romannumeral#1}}
\newcommand{\I}{\text{\Romannumeral{1}}}
\newcommand{\II}{\text{\Romannumeral{2}}}

\newtheorem{theo}{\bf Theorem}[section]

\newcommand{\ma}{\hspace{0pt}}

\allowdisplaybreaks[4]

\begin{document}

\baselineskip=20pt
\newcommand{\preprint}{
\vspace*{-20mm}\begin{flushleft}\end{flushleft}
}
\newcommand{\Title}[1]{{\baselineskip=26pt
  \begin{center} \Large \bf #1 \\ \ \\ \end{center}}}
\newcommand{\Author}{\begin{center}
  \large \bf 
  Ryu Sasaki${}$ \end{center}}
\newcommand{\Address}{\begin{center}
     Department of Physics, Tokyo University of Science,
     Noda 278-8510, Japan
        \end{center}}
\newcommand{\Accepted}[1]{\begin{center}
  {\large \sf #1}\\ \vspace{1mm}{\small \sf Accepted for Publication}
  \end{center}}

\preprint
\thispagestyle{empty}

\Title{Exactly solvable  discrete time Birth and Death processes
}

\Author

\Address
\vspace{1cm}

\begin{abstract}
We present 15 explicit examples of discrete time Birth and Death processes which are {\em exactly solvable}.
They are related to the  hypergeometric orthogonal polynomials of  Askey scheme having discrete orthogonality measures.
Namely, they are the Krawtchouk, three different kinds of $q$-Krawtchouk, (dual, $q$)-Hahn, ($q$)-Racah,
Al-Salam-Carlitz $\II$, $q$-Meixner, $q$-Charlier, dual big $q$-Jacobi and dual big $q$-Laguerre polynomials.
The birth and death rates are determined by the difference equations governing the polynomials.
The stationary distributions are the normalised orthogonality measures of the polynomials.
The transition probabilities are neatly expressed by the normalised polynomials and the corresponding eigenvalues.
This paper is simply the discrete time versions of the known solutions of the continuous time birth and death processes.
\end{abstract}

%
%
\section{Introduction}
It is known \cite{bdsol}, \cite{os34}\S6 that all hypergeometric orthogonal polynomials of a discrete variable
\cite{nikiforov}
belonging to Askey scheme \cite{askey}--\cite{gasper} provide exactly solvable continuous time 
Birth and Death (BD) processes in one dimension.
A good part of these polynomials also supply exactly solvable {\em discrete time} BD processes as will be shown in this paper.
The selection criterion is the boundedness of the birth plus death rate $B(x)+D(x)<\infty$, in which 
birth rate $B(x)$ and death rate $D(x)$ are identified as the coefficients of the difference equations governing these polynomials.
 
 The birth and death processes, continuous and discrete time, are typical examples of stationary 
 Markov processes and chains with
 a wide range of applications \cite{feller}; demography, queueing theory, inventory models, infections
  and chemical dynamics, etc.
 In this paper, however,   only the mathematical sides of the BD processes, the problem setting and solution procedures,
 etc are expanded.  The main point of the logic is that the matrices of transition probabilities of continuous 
 and discrete time BD, the difference equations governing the polynomial, a real symmetric non-negative tri-diagonal 
 matrix equation determining the polynomial \cite{os12} are all connected by similarity transformations in
 terms of a diagonal matrix determined by the {\em stationary probability distribution}, {\em i.e.} the
 orthogonality measure of the polynomial. 
 In most literature on BD processes \cite{ismail,KarMcG,schoutens,bdp}, however,
 the implementation of orthogonal polynomials depended on the three term recurrence relations 
 rather than the difference equations.
 This is why the present simple solution method have not been noticed so long time.
 
 The present paper is prepared in a plain style so that non-experts can easily understand.
 It is  organised as follows.
 An elementary introduction of Markov chains and the problem setting of discrete time BD is given in section two.
 The solution procedures of discrete time BD  are developed in parallel with those of continuous time BD in section three.
 This is the main part of the paper. 
 To an arbitrary continuous time BD with bounded $B(x)+D(x)<\infty$, a discrete time BD is associated with a
 free parameter $t_S$ representing the time scale. 
 As shown in {\bf Theorem 3.1} they share common eigenvectors and the eigenvalues 
 are linearly related \eqref{dctrel},\eqref{eigrel}. 
 A tri-diagonal matrix $\widetilde{\mathcal H}$ \eqref{htdef} is introduced by $B(x)$ and $D(x)$. 
 A positive diagonal matrix $\Phi$ is constructed by the ratios of $B(x)$ and $D(x+1)$, \eqref{phi0def},\eqref{Phidef}.
 A similarity transformation of $\widetilde{\mathcal H}$ in terms of $\Phi$  produces a 
 {\em symmetric and positive-semidefinite tri-diagonal} matrix 
 $\mathcal{H}=\Phi\widetilde{\mathcal H}\Phi^{-1}$ \eqref{Hdef}.
 The complete set of eigenvalues and eigenvectors of $\mathcal{H}$ provide the complete solutions of the
 continuous and discrete time BD, as their transition probability matrices $L_{BD}$ and $L$ are
 also obtained from $\mathcal{H}$ by a similarity transformation, in the opposite direction.
 The solutions of the initial value problem and the transition matrix after time $t$ (step $\ell$)
 are provided in {\bf Theorem 3.2, 3.3}. 
 The spectral representations of the transition matrices are presented in {\bf Theorem 3.4}.
 When the birth and death rates $B(x)$ and $D(x)$ are chosen as the coefficients of the difference equations 
 governing the polynomial \eqref{bdeq}, the corresponding BD processes, continuous and discrete, are
 exactly solvable as demonstrated in {\bf Theorem 3.5, 3.6, 3.7}.
 Various data for exactly solvable BD's, $B(x)$, $D(x)$, the eigenvalues, eigenvectors, etc are 
 presented in section four and five.
 Section five deals with exactly solvable semi-infinite cases in which the complete eigenvectors consist of 
 two sets of mutually orthogonal polynomials.
 Section six discusses two exactly solvable finite cases, which are mirror symmetric at the mid point.
 Section seven is for comments.
%
%
\section{Markov chain}
\label{sec:Markov}
The subject of the present paper belongs to the simplest category of stationary Markov chain
on a one-dimensional integer lattice $\mathcal{X}$, either finite or semi-infinite:
\begin{equation*}
\mathcal{X}=\{0,1,\ldots,N\}: \quad \text{finite},\qquad \mathcal{X}
=\mathbb{Z}_{\ge0}: \quad \text {semi-infinite}.
\end{equation*}
For analytic treatment we use $x,y,..$ as representing the lattice points in 
$\mathcal{X}$, $x,y\in\mathcal{X}$.
The general problem setting is as follows.
Suppose a non-negative matrix $L$  of transition probability is given. Its element
\begin{equation}
L_{x\,y}\ge0,\qquad \sum_{x\in\mathcal{X}}L_{x\,y}=1,
\label{ldef}
\end{equation}
is the transition probability {\em from $y$ to $x$}. 
When the system has the probability distribution
\begin{equation*}
\mathcal{P}(x;\ell)\ge0,\quad \sum_{x\in\mathcal{X}}\mathcal{P}(x;\ell)=1,
\end{equation*}
at $\ell$-th step, the next step distribution is given by
\begin{equation}
\mathcal{P}(x;\ell+1)=\sum_{y\in\mathcal{X}}L_{x\,y}\mathcal{P}(y;\ell).
\label{nextstep}
\end{equation}
The above condition \eqref{ldef}  ensures the conservation of probability,
\begin{equation*}
\sum_{x\in\mathcal{X}}\mathcal{P}(x;\ell+1)=\sum_{y\in\mathcal{X}}
\sum_{x\in\mathcal{X}}L_{x\,y}\mathcal{P}(y;\ell)=\sum_{y\in\mathcal{X}}\mathcal{P}(y;\ell)=1.
\end{equation*}
Another immediate consequence of \eqref{ldef} and Perron-Frobenius theorem
applied to $L$ is that its  spectrum is bounded by $1$ and $-1$:
\begin{equation}
-1\le\text{Eigenvalues}(L)\le1.
\label{eigbound}
\end{equation}
This can be easily seen by considering the normalised eigenvector $v_M$ of $L$ corresponding 
to the maximal eigenvalue $\kappa_M$,
\begin{equation*}
\sum_{y\in\mathcal{X}}L_{x\,y}v_M(y)=\kappa_Mv_M(x),\quad \sum_{y\in\mathcal{X}}|v_M(y)|=1.
\end{equation*}
 Since we can always choose $v_M(y)\ge0$ by Perron-Frobenius theorem, $v_M(y)$ is a probability distribution,
 $\sum_{y\in\mathcal{X}}v_M(y)=1$.
This means that $\kappa_M=1$ since $\kappa_Mv_M(x)$ is also a probability distribution, as shown above.
The lower bound $-1$ in \eqref{eigbound} is derived by applying Perron-Frobenius theorem to $L^2$.
Other eigenvectors of $L$, having sign changes, can never constitute a probability distribution on their own.

Among many problems, the following three are most basic.\\
$\bullet$ {\bf Initial value problem}. Given an initial probability distribution 
\begin{equation*}
\mathcal{P}(x;0)\ge0,\quad \sum_{x\in\mathcal{X}}\mathcal{P}(x;0)=1,
\end{equation*}
calculate the distribution after $\ell$ steps,
\begin{equation}
\mathcal{P}(x;\ell)=\sum_{y\in\mathcal{X}}(L^\ell)_{x\,y}\mathcal{P}(y;0).
\label{inipro}
\end{equation}
$\bullet$ {\bf $\ell$ step transition probability}. 
Starting from the initial distribution concentrated at $y$, $\mathcal{P}(x;0)=\delta_{x\,y}$, 
derive the explicit form of $\ell$ step transition probability from $y$ to $x$
\begin{equation}
\mathcal{P}(x,y;\ell)=(L^\ell)_{x\,y}.
\label{ellpro}
\end{equation}
$\bullet$ {\bf Spectral representation of $L$ in terms of the eigenvalues and the eigenvectors}.

\bigskip
Hereafter let us restrict the matrix $L$ to be {\em tri-diagonal},
\begin{equation}
L_{x\,y}=0 \quad \text{if}\ |x-y|\ge2,
\label{jacobi}
\end{equation}
a {\em non-negative Jacobi matrix}. 
The tri-diagonal restriction means that, at each step, transitions are restricted within 
 the nearest neighbour lattice points. Such Markov chains are usually called generalised random walks. 
 Let us specify the tri-diagonal transition matrix $L$ by using the language of random walk.
 This is achieved by choosing two positive functions of $x$. A walker at point $x$ advances to $x+1$ 
 with  probability $\bar{B}(x)$  and he retreats to $x-1$ with probability $\bar{D}(x)$ and stays at $x$
 with probability $1-\bar{B}(x)-\bar{D}(x)$, {\em i.e.}
 \begin{align}
L_{x+1\,x}=\bar{B}(x),\ \  L_{x-1\,x}=\bar{D}(x),\ \ L_{x\,x}=1-\bar{B}(x)-\bar{D}(x), 
\label{Lndef1}\\
 0<\bar{B}(x)<1,\quad 0< \bar{D}(x)<1,\quad 0< \bar{B}(x)+\bar{D}(x)<1,
\end{align}
 together  with the boundary condition(s) 
\begin{equation}
\bar{D}(0)=0,\quad  \bar{B}(N)=0:\  (\text{only  for  finite cases}).
\label{bdcond1}
\end{equation} 
The matrix $L$ looks as follows:
\begin{align*}
&L=\n
&{\small
\left(
\begin{array}{cccccc}
1-\bar{B}(0)  & \bar{D}(1)  &   0& \cdots&\cdots&0\\
\bar{B}(0)  & 1-\bar{B}(1)-\bar{D}(1)  &\bar{D}(2) &0&\cdots&\vdots  \\
\!\!0  &  \bar{B}(1)  &  1-  \bar{B}(2)-\bar{D}(2)&\bar{D}(3)&\cdots&\vdots\\
\!\!\vdots&\cdots&\cdots&\cdots&\cdots&\vdots\\
\!\!\vdots&\cdots&\cdots&\cdots&\cdots&0\\
\!\!0&\cdots&\cdots&\bar{B}(N\!\!-\!\!2)&1-\bar{B}(N\!\!-\!\!1)\!-\bar{D}(N\!\!-\!\!1)&\bar{D}(N)\\
\!\!0&\cdots&\cdots&0&\bar{B}(N\!-\!1)&1-\bar{D}(N)
\end{array}
\right).
}
\end{align*}
The equation \eqref{nextstep} connecting the $\ell$-th step distribution and $\ell+1$-th step distribution
reads for the inner points $1\le x\le N-1$, 
\begin{equation}
\mathcal{P}(x;\ell+1)=
(1-\bar{B}(x)-\bar{D}(x))\mathcal{P}(x;\ell)+\bar{B}(x-1)\mathcal{P}(x-1;\ell)+\bar{D}(x+1)\mathcal{P}(x+1;\ell),
\end{equation}
and for the endpoint(s)
\begin{align}
\mathcal{P}(0;\ell+1)&=(1-\bar{B}(0))\mathcal{P}(0;\ell)+\bar{D}(1)\mathcal{P}(1;\ell),\\
\mathcal{P}(N;\ell+1)&=\bar{B}(N-1)\mathcal{P}(N-1;\ell)+
(1-\bar{D}(N))\mathcal{P}(N;\ell).
\label{randomend}
\end{align}

%
%
\section{Discrete time  Birth and Death process}
\label{sec:dtbd}
Probably it is now clear that the above generalised random walk \eqref{Lndef1}--\eqref{randomend} 
can also be called {\em discrete time} Birth and Death (BD) process.
They are very closely related to the {\em continuous time} Birth and Death process.
For comparison, let us review it here.
Let $\mathcal{P}(x;t)$ be the probability distribution over $\mathcal{X}$ at time $t$.
Let us denote the birth rate at population $x$ by $B(x)>0$ and the death rate by $D(x)>0$.
The time evolution of the probability distribution is governed by the following differential equation:
\begin{align}
\frac{\partial}{\partial t}\mathcal{P}(x;t)&=(L_{BD}\mathcal{P})(x;t)=\sum_{y\in\mathcal{X}}{L_{BD}}_{x\,y}\mathcal{P}(y;t),
\quad \mathcal{P}(x;t)\ge0,\quad \sum_{x\in\mathcal{X}} \mathcal{P}(x;t)=1,
\label{bdeqformal}\\
&=-(B(x)+D(x))\mathcal{P}(x;t)+B(x-1)\mathcal{P}(x-1;t)+D(x+1)\mathcal{P}(x+1;t),
\label{BDeq}
\end{align}
with the boundary condition(s)
\begin{equation}
D(0)=0,\quad B(N)=0:\  (\text{only  for a finite case}),
\label{bdcond2}
\end{equation}
which is called reflecting boundary condition.
Here the matrix $L_{BD}$ is also tri-diagonal
\begin{align}
&{L_{BD}}_{x+1\,x}=B(x),\ \  {L_{BD}}_{x-1\,x}=D(x),\ {L_{BD}}_{x\,x}=-B(x)-D(x), \n
&\hspace{30mm}{L_{BD}}_{x\,y}=0,\quad |x-y|\ge2,
\label{LBDdef}\\
&L_{BD}=\n
&\left(
\begin{array}{cccccc}
-B(0)  & D(1)  &   0& \cdots&\cdots&0\\
B(0)  &-B(1)-D(1)  &D(2) &0&\cdots&\vdots  \\
\!\!0  &  B(1)  &  -B(2)-D(2)&D(3)&\cdots&\vdots\\
\!\!\vdots&\cdots&\cdots&\cdots&\cdots&\vdots\\
\!\!\vdots&\cdots&\cdots&\cdots&\cdots&0\\
\!\!0&\cdots&\cdots&B(N\!\!-\!\!2)&-B(N\!\!-\!\!1)\!-\!D(N\!\!-\!\!1)&D(N)\\
\!\!0&\cdots&\cdots&0&B(N\!-\!1)&-D(N)
\end{array}
\right),\nonumber
\end{align}
satisfying the condition
\begin{equation}
\sum_{x\in\mathcal{X}}{L_{BD}}_{x\,y}=0.
\label{bdproconv}
\end{equation}
This ensures the conservation of probability, that is, the condition
$\sum_{x\in\mathcal{X}} \mathcal{P}(x;t)=1$  is preserved by the time evolution \eqref{bdeqformal}.

\bigskip
By introducing a parameter $t_S$ specifying the time spacing, a discrete time BD process is obtained from a
continuous time BD process as shown in the following
\begin{theo}
\label{theo1}
For each continuous BD process with {\em bounded} $B(x)+D(x)$, a discrete time BD process is defined
with one free positive parameter $t_S$ as follows,
\begin{align} 
&\bar{B}(x)\eqdef t_SB(x),\quad \bar{D}(x)\eqdef t_SD(x),
\quad t_S\cdot{\rm max}\left(B(x)+D(x)\right)<1,
\label{dctrel}\\
&\Longrightarrow L=I_d+t_S\cdot L_{BD},\qquad I_d:\ \text{Identity matrix}.
\label{eigrel}
\end{align}
If the continuous time BD is solved, the corresponding discrete time BD is solved, and vice versa, 
as  the eigenvalues are related as
\begin{equation*}
\text{Eigenvalue of}\ L=1+t_S\cdot\left(\text{Eigenvalue of}\ (L_{BD})\right),
\end{equation*}
and the corresponding eigenvectors are common. 
\end{theo}
The bigger $t_S$, the bigger is the time interval of the corresponding discrete time BD process.
Obviously, $t_S$ has an upper limit given by \eqref{dctrel}.
Later it will be shown that the spectrum of $L_{BD}$ is {\em negative semi-definite} \eqref{LBDsol}.

\bigskip
Let us proceed to solve the continuous time BD 
\eqref{bdeqformal}--\eqref{bdcond2} 
in the general setting, {\em i.e.} $B(x)$ and $D(x)$ are arbitrary positive functions restricted only by the 
boundary condition(s) \eqref{bdcond2}. The special cases related to the hypergeometric 
orthogonal polynomials \cite{bdsol} will be discussed later.
Let us introduce a {\em tri-diagonal matrix} $\widetilde{\mathcal H}$ on 
$\mathcal{X}$ in terms of the birth and death rates
$B(x)$ and $D(x)$,
\begin{align} 
 &\widetilde{\mathcal H}_{x\,x+1}=-B(x),\ \  \widetilde{\mathcal H}_{x\,x-1}=-D(x),\ 
\widetilde{\mathcal H}_{x\,x}=B(x)+D(x), \ 
\widetilde{\mathcal H}_{x\,y}=0, |x-y|\ge2,
\label{htdef}\\
&\widetilde{\mathcal H}=\n
&\left(
\begin{array}{cccccc}
\!\!B(0)  & -B(0)  &   0& \cdots&\cdots&0\\
\!\!-D(1)  & B(1)+D(1)  & -B(1) &\cdots&\cdots&\vdots  \\
\!\!0  &  -D(2) &   B(2)+D(2)&\cdots&\cdots&\vdots\\
\!\!\vdots&\cdots&\cdots&\cdots&\cdots&\vdots\\
\!\!\vdots&\cdots&\cdots&\cdots&\cdots&0\\
\!\!0&\cdots&\cdots&\cdots&B(N\!\!-\!\!1)\!+\!D(N\!\!-\!\!1)&-B(N\!\!-\!\!1)\\
\!\!0&\cdots&\cdots&\cdots&-D(N)&D(N)
\end{array}
\right),\nonumber
\end{align}
and consider its eigenvalue problem
\begin{equation}
 (\widetilde{\mathcal H}\check{P}_n)(x)=\sum_{y\in\mathcal{X}}\widetilde{\mathcal H}_{x\,y}\check{P}_n(y)=\mathcal{E}(n)\check{P}_n(x),\ x\in\mathcal{X},\ n\in\mathcal{X}.
  \label{hteq}
\end{equation}
The equation reads explicitly as
\begin{equation}
B(x)\left(\check{P}_n(x)-\check{P}_n(x+1)\right)+
D(x)\left(\check{P}_n(x)-\check{P}_n(x-1)\right)=\mathcal{E}(n)\check{P}_n(x),\quad n\in\mathcal{X}.
\label{bdeq0}
\end{equation}
Here $\check{P}_n(x)$ and $\mathcal{E}(n)$ are to be determined as an eigenvector and the 
corresponding  eigenvalue of $\widetilde{\mathcal H}$.
It is well known that the top component of an eigenvector of a tri-diagonal matrix is non-vanishing.
We adopt the following  universal normalisation of the eigenvectors $\{\check{P}_n(x)\}$,
\begin{equation}
\check{P}_n(0)=1,\quad n\in\mathcal{X}.
\label{pnorm}
\end{equation}
Then  it is obvious
\begin{equation}
\sum_{y\in\mathcal{X}}\widetilde{\mathcal H}_{x\,y}=0, \
\check{P}_0(x)\eqdef1,\  (\forall x\in\mathcal{X}),\quad \Longrightarrow 
\sum_{y\in\mathcal{X}}\widetilde{\mathcal H}_{x\,y}\check{P}_0(y)=0,
\label{constsol}
\end{equation}
namely,  a constant vector of identical components is the eigenvector of $\widetilde{\mathcal H}$ 
of vanishing eigenvalue $\mathcal{E}(0)=0$. 
This is also obvious from the difference equation \eqref{bdeq0}.
The matrix $\widetilde{\mathcal H}$  is related to a {\em real symmetric  
tri-diagonal matrix} $\mathcal{H}$ by a similarity transformation. 
Let us introduce a positive function 
$\phi_0(x)$ on $\mathcal{X}$ and a {\em diagonal matrix} $\Phi$ consisting of $\phi_0(x)$ 
by the ratios of $B(x)$ and $D(x+1)$,
\begin{align} 
  &\phi_0(0)\eqdef1,\quad \phi_0(x)\eqdef\sqrt{\prod_{y=0}^{x-1}\frac{B(y)}{D(y+1)}}
  \Leftrightarrow  \frac{\phi_0(x+1)}{\phi_0(x)}=\frac{\sqrt{B(x)}}{\sqrt{D(x+1)}},\quad x\in\mathcal{X}.
  \label{phi0def}\\
& \Phi_{x\,x}=\phi_0(x),\quad \Phi_{x\,y}=0, \quad x\neq y.
\label{Phidef}
\end{align}
For semi-infinite cases, $B(x)$ and $D(x)$ must be restricted so that $\phi_0(x)$ is square summable,
\begin{equation}
\sum_{x\in\mathcal{X}}\phi_0(x)^2<\infty.
\end{equation}

Let us define
\begin{equation}
\phi_n(x)\eqdef \phi_0(x)\check{P}_n(x),\quad n\in\mathcal{X},
\label{phindef}
\end{equation}
which constitutes the eigenvector of the real symmetric matrix $\mathcal{H}$ defines as follows,
\begin{align} 
   &\mathcal{H}\eqdef \Phi \widetilde{\mathcal H} \Phi^{-1}\ \Leftrightarrow 
   \widetilde{\mathcal H}=\Phi^{-1} \mathcal{H} \Phi \ \Leftrightarrow
   \mathcal{H}_{x\,y}=\phi_0(x)\widetilde{\mathcal H}_{x\,y}\phi_0(y)^{-1},
   \label{Hdef}\\
  &\hspace{60mm} \Downarrow\n
&(\mathcal{H}\phi_n)(x) =\sum_{y\in\mathcal{X}}\mathcal{H}_{x\,y}\phi_n(y)=
\sum_{y\in\mathcal{X}}\phi_0(x)\widetilde{\mathcal H}_{x\,y}\check{P}_n(y)
=\mathcal{E}(n)\phi_n(x),\quad n\in\mathcal{X},
\label{phineq}\\
&{\mathcal H}_{x\,x+1}=-\sqrt{B(x)D(x+1)},\ \  {\mathcal H}_{x\,x-1}=-\sqrt{B(x-1)D(x)},\ 
{\mathcal H}_{x\,x}=B(x)+D(x), \n
&\hspace{70mm}{\mathcal H}_{x\,y}=0, \quad |x-y|\ge2,
\label{tdefcomp}\\
&\mathcal{H}=\n
&\left(
\begin{array}{cccccc}
\!\!B(0)  & -\sqrt{B(0)D(1)}  &   0& \cdots&\cdots&0\\
\!\!-\sqrt{B(0)D(1)}  & B(1)+D(1)  & -\sqrt{B(1)D(2)} &\cdots&\cdots&\vdots  \\
\!\!0  &  -\sqrt{B(1)D(2)}  &   B(2)+D(2)&\cdots&\cdots&\vdots\\
\!\!\vdots&\cdots&\cdots&\cdots&\cdots&\vdots\\
\!\!\vdots&\cdots&\cdots&\cdots&\cdots&0\\
\!\!0&\cdots&\cdots&\cdots&B(N\!\!-\!\!1)\!+\!D(N\!\!-\!\!1)&-\sqrt{\!B(N\!\!-\!\!1)D(N)}\\
\!\!0&\cdots&\cdots&0&-\sqrt{B(N\!-\!1)D(N)}&D(N)
\end{array}
\right).
\end{align}
The real symmetry and the positive semi-definiteness of $\mathcal{H}$ can be seen clearly 
by the following factorisation in terms of an upper triangular matrix $\mathcal{A}$,
\begin{equation}
\mathcal{H}={}^t\!\!\mathcal{A}\mathcal{A}\ \Rightarrow \mathcal{H}
={}^t\mathcal{H},\quad \mathcal{A}_{x\,x}=\sqrt{B(x)},\quad 
\mathcal{A}_{x\,x+1}=-\sqrt{D(x+1)},\quad \mathcal{A}_{x\,y}=0,\ \text{otherwise},
\end{equation}
in which ${}^t\!\!\mathcal{A}$ is the transposed matrix of $\mathcal{A}$.
This guarantees the reality and non-negativeness of the eigenvalues 
\begin{equation}
\mathcal{E}(n)\ge0,\qquad n\in\mathcal{X},
\label{posE}
\end{equation}
and the orthogonality of the eigenvectors $\{\phi_n(x)\}$ of $\mathcal{H}$, since the  simpleness of
the eigenvalues is due to its tri-diagonality. 
It should be stressed that $\phi_0(x)$ is the zero mode (eigenvector) of $\mathcal{A}$ and $\mathcal{H}$,
\begin{equation}
0=(\mathcal{A}\phi_0)(x)=\sqrt{B(x)}\phi_0(x)-\sqrt{D(x+1)}\phi_0(x+1)
\ \Rightarrow (\mathcal{H}\phi_0)(x)=0,
\label{Azero}
\end{equation}
and $\phi_0(x)^2$ provides the orthogonality measure of the eigenvectors  $\{\check{P}_n(x)\}$ 
of $\widetilde{\mathcal H}$,
\begin{align} 
(\phi_n,\phi_m)\eqdef \sum_{x\in\mathcal{X}}\phi_n(x)\phi_m(x)
= \sum_{x\in\mathcal{X}}\phi_0(x)^2\check{P}_n(x)\check{P}_m(x)
=\frac1{d_n^2}\delta_{n\,m},
 \quad n, m\in\mathcal{X},
  \label{orth}
\end{align}
in which the normalisation constants $\{d_n>0\}$ are calculated 
after all the eigenvectors $\{\check{P}_n(x)\}$ are known. 
It should be stressed that $\phi_0(x)$ and $\check{P}_n(x)$ are uniquely specified by the
normalisation condition \eqref{pnorm}, \eqref{phi0def}
\begin{equation*}
\phi_0(0)=1=\check{P}_n(0),\quad n\in\mathcal{X}.
\end{equation*}

Let us define orthonormal vectors $\{\hat{\phi}_n(x)\}$ 
\begin{equation}
\hat{\phi}_n(x)\eqdef d_n\phi_n(x)=d_n\phi_0(x)\check{P}_n(x),\quad (\hat{\phi}_n,\hat{\phi}_m)=\delta_{n\,m},
\quad n,m \in\mathcal{X},
\label{normphin}
\end{equation}
and the square of the normalised zero mode
\begin{align}
\pi(x)&\eqdef \hat{\phi}_0(x)^2=d_0^2\phi_0(x)^2=d_0^2\,\prod_{y=0}^{x-1}\frac{B(y)}{D(y+1)},
\label{pidef}\\
& \sum_{x\in\mathcal{X}}\pi(x)=1 \Leftarrow \frac1{d_0^2}\eqdef\sum_{x\in\mathcal{X}}\prod_{y=0}^{x-1}\frac{B(y)}{D(y+1)}=\sum_{x\in\mathcal{X}}\prod_{y=0}^{x-1}\frac{\bar{B}(y)}{\bar{D}(y+1)},
\end{align}
which will turn out to be the {\em stationary distribution}.

\bigskip
With these preparations, let us return to $L_{BD}$ \eqref{LBDdef}.
It is now clear that $L_{BD}$ is also related to $\mathcal{H}$ by a  similarity transformation 
in terms of $\Phi$, but in the opposite direction to $\widetilde{\mathcal H}$ together with a negative sign,
\begin{align} 
   &\hspace{50mm}L_{BD}=-\Phi\mathcal{H}\Phi^{-1},
   \label{LBDHrel} \\
& {L_{BD}}_{x\,x}=-\mathcal{H}_{x\,x}=-B(x)-D(x),\n
&{L_{BD}}_{x+1\,x}=-\phi_0(x+1)\mathcal{H}_{x+1\,x}\phi_0(x)^{-1}
=\phi_0(x+1)\sqrt{B(x)D(x+1)}\,\phi_0(x)^{-1}
=B(x),\n
&{L_{BD}}_{x-1\,x}=-\phi_0(x-1)\mathcal{H}_{x-1\,x}\phi_0(x)^{-1}
=\phi_0(x-1)\sqrt{B(x-1)D(x)}\,\phi_0(x)^{-1}
=D(x),\n[4pt]
&\hspace{40mm} (L_{BD}\hat{\phi}_0\hat{\phi}_n)(x)=-\mathcal{E}(n)\hat{\phi}_0(x)\hat{\phi}_n(x),\quad
n\in\mathcal{X}.
\label{LBDsol}
\end{align}
Thus we arrive at the solutions of the general continuous time BD by the following
\begin{theo}
\label{theo3}
If we obtain the complete set of eigensystem of the matrix $\widetilde{\mathcal H}$ \eqref{htdef}--\eqref{hteq}, 
the solution of the initial value problem of the continuous time BD  \eqref{bdeqformal}--\eqref{LBDdef} 
is given by
\begin{equation}
\mathcal{P}(x;t)=\hat{\phi}_0(x)\sum_{n\in\mathcal{X}}c_ne^{-\mathcal{E}(n)t}\hat{\phi}_n(x),
\label{ctbdsol1}
\end{equation}
in which $\{c_n\}$ are determined as the expansion coefficients of the initial distribution $\mathcal{P}(x;0)$,
\begin{equation}
\mathcal{P}(x;0)=\hat{\phi}_0(x)\sum_{n\in\mathcal{X}}c_n\hat{\phi}_n(x)\ \Rightarrow
c_0=1,\quad c_n=\sum_{x\in\mathcal{X}}\hat{\phi}_n(x)\hat{\phi}_0(x)^{-1}\mathcal{P}(x;0),
\quad n=1,\ldots.
\label{cndef}
\end{equation}
The transition matrix from $y$ to $x$ after time $t$ is
\begin{equation}
\mathcal{P}(x,y;t)=\hat{\phi}_0(x)\hat{\phi}_0(y)^{-1}
\sum_{n\in\mathcal{X}}e^{-\mathcal{E}(n)t}\hat{\phi}_n(x)\hat{\phi}_n(y).
\end{equation}
The approach to the stationary distribution $\pi(x)$ \eqref{pidef}  is guaranteed by the positivity of $\mathcal{E}(n)>0$
for the non-zero modes,
\begin{equation}
\lim_{t\to\infty}\mathcal{P}(x;t)=\pi(x),\quad \lim_{t\to\infty}\mathcal{P}(x,y;t)=\pi(x).
\label{piapp}
\end{equation}
\end{theo}
Based on the relationship between the discrete and continuous BD's {\bf Theorem \ref{theo1}} \eqref{dctrel}--\eqref{eigrel}
we arrive at the solutions of the general discrete time BD by the following
\begin{theo}
\label{theo4}
The complete set of eigensystem of the matrix $\widetilde{\mathcal H}$ \eqref{hteq}
provides the complete set of eigensystem of the matrix $L$ \eqref{Lndef1} for the discrete time BD 
with the identification \eqref{dctrel}
\begin{equation}
(L\hat{\phi}_0\hat{\phi}_n)(x)=\kappa(n)\hat{\phi}_0(x)\hat{\phi}_n(x),
\quad \kappa(n)=1-{t_S}\cdot\mathcal{E}(n),\quad n\in\mathcal{X}.
\end{equation}
The  solution of the initial value problem of the discrete time BD  \eqref{Lndef1}--\eqref{bdcond1} 
after $\ell$ steps is given by
\begin{equation}
\mathcal{P}(x;\ell)=\hat{\phi}_0(x)\sum_{n\in\mathcal{X}}c_n\kappa(n)^\ell\hat{\phi}_n(x),
\label{dtbdsol1}
\end{equation}
in which $\{c_n\}$ are given in \eqref{cndef}. The $\ell$ step transition matrix from $y$ to $x$ is
\begin{equation}
\mathcal{P}(x,y;\ell)=\hat{\phi}_0(x)\hat{\phi}_0(y)^{-1}
\sum_{n\in\mathcal{X}}\kappa(n)^\ell\hat{\phi}_n(x)\hat{\phi}_n(y).
\label{dtbdsol2}
\end{equation}
The approach to the stationary distribution is about the same as \eqref{piapp}.
It should be stressed that for both continuous and discrete time BD's, the stationary distribution $\pi(x)$ 
is the same
and it is determined by the input $B(x), D(x)$ and $\bar{B}(x),\bar{D}(x)$ only without solving the
eigenvalue problem of $\widetilde{\mathcal H}$ \eqref{hteq}.
\end{theo}
\begin{theo}
\label{theoSpe}
The complete set of eigensystem of the matrix $\widetilde{\mathcal H}$ \eqref{hteq}
provides the spectral representation of the symmetric matrix $\mathcal{H}$ \eqref{Hdef},
\begin{equation}
\mathcal{H}_{x\,y}=\sum_{n\in\mathcal{X}}\mathcal{E}(n)\hat{\phi}_n(x)\hat{\phi}_n(y).
\label{Hspe}
\end{equation}
This in turn supplies the spectral representations of $L_{BD}$ and $L$ through \eqref{LBDHrel} and \eqref{eigrel},
\begin{align} 
 {L_{BD}}_{x\,y}&=-\hat{\phi}_0(x)\sum_{n\in\mathcal{X}}\mathcal{E}(n)
 \hat{\phi}_n(x)\hat{\phi}_n(y)\hat{\phi}_0(y)^{-1}
 =-\pi(x)\sum_{n\in\mathcal{X}}\mathcal{E}(n)({d_n^2}/{d_0^2})\check{P}_n(x)\check{P}_n(y),\\
 L_{x\,y}&=\hat{\phi}_0(x)\sum_{n\in\mathcal{X}}\kappa(n)\hat{\phi}_n(x)\hat{\phi}_n(y)\hat{\phi}_0(y)^{-1}
=\pi(x)\sum_{n\in\mathcal{X}}\kappa(n)({d_n^2}/{d_0^2})\check{P}_n(x)\check{P}_n(y),
\end{align}
in which $\pi(x)$ \eqref{pidef} is the stationary distribution.
\end{theo}

As for the solutions of  the above {\em continuous time} Birth and Death process 
\eqref{bdeqformal}--\eqref{bdproconv}, we have the following 
\begin{theo}
\label{theo2}
The above continuous time Birth and Death process \eqref{bdeqformal}--\eqref{LBDdef} is 
{\bf exactly solvable} when $B(x)$ and $D(x)$ are chosen to be the
coefficient functions of certain difference equations \cite{bdsol}
\begin{align}
&B(x)\left(\check{P}_n(x)-\check{P}_n(x+1)\right)+
D(x)\left(\check{P}_n(x)-\check{P}_n(x-1)\right)=\mathcal{E}(n)\check{P}_n(x),\quad n\in\mathcal{X},
\label{bdeq}\\
&\hspace{10mm} \Longleftrightarrow \quad \sum_{y\in\mathcal{X}}\widetilde{\mathcal H}_{x,y}
\check{P}_n(y)=\mathcal{E}(n)\check{P}_n(x),\quad n\in\mathcal{X},\nonumber
\end{align}
which determine  hypergeometric orthogonal polynomials
$\{\check{P}_n(x)\}$ with discrete orthogonality measures belonging to Askey scheme.
In this case all the eigenvalues of $L_{BD}$ and the corresponding eigenvectors are explicitly known.
The eigenvalues of $L_{BD}$ are the same as those in \eqref{bdeq} 
with a minus sign $\{-\mathcal{E}(n)\}$ $n=0,1\ldots$,
and the eigenvectors are proportional to $\{\check{P}_n(x)\}$ of \eqref{bdeq}, as explicitly given in \eqref{LBDsol}.
$\{\check{P}_n(x)\}$ are hypergeometric orthogonal
polynomial 
\begin{equation}
\check{P}_n(x)=P_n\bigl(\eta(x)\bigr),\quad n\in\mathcal{X},
\label{pnsin}
\end{equation}
in a certain {\em sinusoidal coordinate} $\eta(x)$, which takes the following five types \cite{os12}.
They are linear  or quadratic in $x$ and linear or `quadratic' in $q^{\pm x}$, with $0<q<1$,
\begin{equation}
\eta(x): \quad x, \quad x(x+d),\quad 1-q^x,\quad q^{-x}-1, \quad (q^{-x}-1)(1-dq^x);
\quad \eta(0)=0.
\end{equation}
It should be stressed that, except for the cases of $\eta(x)=x$, $\check{P}_n(x)$ is not a degree $n$ polynomial in $x$. 
\end{theo}
\begin{theo}
\label{theo22}
The discrete time BD \eqref{Lndef1}--\eqref{bdcond1} is {\bf exactly solvable} if $\bar{B}(x)$ and $\bar{D}(x)$ 
are related to $B(x)$ and $D(x)$ of an exactly solvable continuous time BD by the relation \eqref{dctrel}.
The general forms of the solutions of the initial value problem and the transition matrix will be presented in 
{\bf Theorem \ref{theo5}} for both  continuous and discrete time BD processes.
\end{theo}

\noindent{\bf Remark}
{\bf Theorem \ref{theo2}} was proven in \cite{bdsol} for the continuous time BD processes related to 16 different  polynomials. 
Some of them have unbounded $B(x)+D(x)$.
The   continuous time BD processes related to the polynomials having Jackson integral type measures
are defined on a direct sum of two semi-infinite integer lattices, 
$\mathcal{X}=\mathbb{Z}_{\ge0}\oplus\mathbb{Z}_{\ge0}$, and they require a different formalism.
The big $q$-Jacobi is the typical example.
The solutions of those BD processes are given in \cite{os34}. 
Those polynomials have no corresponding discrete time BD processes as they all have unbounded $B(x)+D(x)$. 
However, the dual polynomials of some of them, {\em e.g.} dual big $q$-Jacobi and dual big $q$-Laguerre
have the above continuous time BD processes \eqref{bdeqformal}--\eqref{LBDdef}, whose solutions are also provided in \cite{os34}.
The solutions of the discrete time BD processes corresponding to  dual big $q$-Jacobi and dual big $q$-Laguerre
are given in \S\ref{dbqj} and \S\ref{dbql}.
The solutions of the  continuous time BD processes related with $q$-Meixner and $q$-Charlier  reported in
\cite{bdsol} were flawed due to the lack of completeness of those polynomials listed in the literature \cite{koeswart,gasper}.
The complete solutions are given in \cite{os34}.
The solutions of the discrete time BD processes corresponding to $q$-Meixner and $q$-Charlier are
listed in \S\ref{qMei} and \S\ref{qcha}. \hfill $\triangle$

\begin{theo}
\label{theo5}
For the exactly solvable continuous and discrete time BD's stated in {\bf Theorem \ref{theo2}, \ref{theo22}},
the formulas of the solutions for the initial value problem, the transition matrix and the spectral representations 
are  the same as those
given in {\bf Theorem \ref{theo3}, \ref{theo4}, \ref{theoSpe}}, so far as the polynomials determined by the difference equations
\eqref{bdeq} form a complete set. 
\end{theo}
{\bf Remark} The examples presented in the next section \S\ref{ex1} belong to this category.
All quantities appearing in the formulas in {\bf Theorem \ref{theo3}, \ref{theo4}, \ref{theoSpe}} 
are explicitly known and reported in the subsections having the names of the corresponding polynomials.
In contrast, the solutions of the exactly solvable BD's presented in section \ref{ex2} require another set of
polynomials on top of those determined by $B(x)$ and $D(x)$ \eqref{bdeq} for completeness.
The extra set of polynomials is also associated with another BD process with its birth and death probabilities
$B^{(-)}(x)$ and $D^{(-)}(x)$, which are related to the original  $B(x)$ and $D(x)$ 
by certain parameter transformations (involutions).
The explicit formulas for these two discrete time BD's are given in {\bf Theorem \ref{theo51}, \ref{theo52},
\ref{theoSpe2}, \ref{theo53}} in the subsections having the names of the polynomials.
The formulas for the corresponding continuous time BD's are not listed since they 
are already reported in \cite{os34} \S6.A.

%
%
\section{Explicit Examples I}
\label{ex1}

Here we present the data for exactly solvable discrete time BD processes, 
for which $B(x)+D(x)$ is bounded and the corresponding polynomials
are complete. 
The ranges of parameters in the birth and death rates are restricted by the positivity of $B(x)$ and $D(x)$.
We list a representative one only. For more general information of the polynomials, 
we refer to \cite{os12} and \cite{koeswart}.
The format for  the normalisation constant $d_n^2$
consists of two parts separated by a $\times$ symbol:
$d_n^2=(d_n^2/d_0^2)\times d_0^2$.  The second part $d_0^2$ satisfies the
relation $\sum_x\phi_0(x)^2=1/d_0^2$.
Throughout sections \ref{ex1} and \ref{ex2}, the parameter $q$ is $0<q<1$.

\subsection{Krawtchouk}
\label{kra}
The case of linear birth and death rates is a very well-known example
(the Ehrenfest model) 
of an exactly solvable birth and death processes \cite{feller, schoutens}:
\begin{align}
  B(x)&=p(N-x),\quad
  D(x)=(1-p)x,\quad 0<p<1,
  \label{kra1}\\
  \mathcal{E}(n)&=n,\qquad
  \eta(x)=x,\\
  \phi_0(x)^2&=
  \frac{N!}{x!\,(N-x)!}\Bigl(\frac{p}{1-p}\Bigr)^x,\quad
  d_n^2
  =\frac{N!}{n!\,(N-n)!}\Bigl(\frac{p}{1-p}\Bigr)^n\times(1-p)^N,\\[4pt]
    \check{P}_n(x)&=P_n(\eta(x))
  ={}_2F_1\Bigl(
  \genfrac{}{}{0pt}{}{-n,\,-x}{-N}\Bigm|p^{-1}\Bigr). 
  \label{kra4}
\end{align}
The stationary probability 
$\pi(x)=\phi_0(x)^2d_0^2=\binom{N}{x}p^x(1-p)^{N-x}$ is the binomial distribution.

\subsection{Hahn}
\label{hah}
This is  a well-known example of quadratic (in $x$) birth and death rates,
\begin{equation}
B(x)=(x+a)(N-x),\quad
  D(x)= x(b+N-x), \quad a>0,\ b>0.
  \label{hah1}
  \end{equation}
It has a quadratic energy spectrum
\begin{align}
\mathcal{E}(n)&= n(n+a+b-1),\quad
  \eta(x)=x,\quad
\phi_0(x)^2
 =\frac{N!}{x!\,(N-x)!}\,\frac{(a)_x\,(b)_{N-x}}{(b)_N},\\
  d_n^2
  &=\frac{N!}{n!\,(N-n)!}\,
  \frac{(a)_n\,(2n+a+b-1)(a+b)_N}{(b)_n\,(n+a+b-1)_{N+1}}
  \times\frac{(b)_N}{(a+b)_N},\\[4pt]
  \check{P}_n(x)&= P_n(\eta(x))
  ={}_3F_2\Bigl(
  \genfrac{}{}{0pt}{}{-n,\,n+a+b-1,\,-x}
  {a,\,-N}\Bigm|1\Bigr).
  \label{hah4}
 \end{align}

\subsection{dual Hahn}
\label{dhah}
The birth and death rates are rational functions of $x$, with $a>0$, $b>0$,
\begin{align}
&B(x)=\frac{(x+a)(x+a+b-1)(N-x)}
  {(2x-1+a+b)(2x+a+b)},
\quad
  D(x)=\frac{x(x+b-1)(x+a+b+N-1)}
  {(2x-2+a+b)(2x-1+a+b)},
  \label{dualhahnBD2}\\
&\mathcal{E}(n)=n,\quad
  \eta(x)= x(x+a+b-1),\quad
  \phi_0(x)^2
  =\frac{N!}{x!\,(N-x)!}
  \frac{(a)_x\,(2x+a+b-1)(a+b)_N}{(b)_x\,(x+a+b-1)_{N+1}},
  \label{dualhahneeta}\\
    &\qquad\qquad d_n^2
  =\frac{N!}{n!\,(N-n)!}\,\frac{(a)_n\,(b)_{N-n}}{(b)_N}
  \times\frac{(b)_{N}}{(a+b)_N},\\[4pt]
    &\check{P}_n(x)=P_n(\eta(x))
  ={}_3F_2\Bigl(
  \genfrac{}{}{0pt}{}{-n,\,x+a+b-1,\,-x}
  {a,\,-N}\Bigm|1\Bigr). 
 \end{align}

\subsection{Racah}
\label{rac}
The function $B(x)$ and $D(x)$ depend on four real parameters $a$, $b$, $c$ and $d$,
with one of them, say $c$,  being related to $N$, $c\equiv -N$:
\begin{align}
&B(x)
  =-\frac{(x+a)(x+b)(x+c)(x+d)}{(2x+d)(2x+1+d)},\quad
  D(x)
  =-\frac{(x+d-a)(x+d-b)(x+d-c)x}{(2x-1+d)(2x+d)},
\label{racahbd}\\
 &\qquad\qquad\qquad\qquad\qquad  a\ge b,\quad  d>0,\quad a>N+d,\quad 0<b<1+d,\\
& \qquad\mathcal{E}(n)= n(n+\tilde{d}),\quad
  \eta(x)=x(x+d),\quad 
  \tilde{d}\eqdef a+b+c-d-1,\\
 & \phi_0(x)^2=\frac{(a,b,c,d)_x}{(1+d-a,1+d-b,1+d-c,1)_x}\,
  \frac{2x+d}{d},\\[4pt]
&d_n^2=\frac{(a,b,c,\tilde{d})_n}
  {(1+\tilde{d}-a,1+\tilde{d}-b,1+\tilde{d}-c,1)_n}\,
  \frac{2n+\tilde{d}}{\tilde{d}}\times
  \frac{(-1)^N(1+d-a,1+d-b,1+d-c)_N}{(\tilde{d}+1)_N(d+1)_{2N}},\\
  &\check{P}_n(x)=P_n(\eta(x))
  ={}_4F_3\Bigl(
  \genfrac{}{}{0pt}{}{-n,\,n+\tilde{d},\,-x,\,x+d}
  {a,\,b,\,c}\Bigm|1\Bigr).\
 \end{align}

\subsection{affine $q$-Krawtchouk}
The birth and death rates are
quadratic  in $q^x$:
\begin{align}
  B(x)&=(q^{x-N}-1)(1-pq^{x+1}),\quad
  D(x)=pq^{x-N}(1-q^x),\quad 0<p<q^{-1}, \\
  \mathcal{E}(n)&=q^{-n}-1,\qquad
  \eta(x)=q^{-x}-1,\\
  \phi_0(x)^2&=\frac{(q\,;q)_N}{(q\,;q)_x(q\,;q)_{N-x}}\,
  \frac{(pq\,;q)_x}{(pq)^x}\,,\quad
  d_n^2
  =\frac{(q\,;q)_N}{(q\,;q)_n(q\,;q)_{N-n}}\,
  \frac{(pq\,;q)_n}{(pq)^n}\times(pq)^N,\\[4pt]
 \check{P}_n(x)&=  P_n(\eta(x))
  ={}_3\phi_2\Bigl(
  \genfrac{}{}{0pt}{}{q^{-n},\,q^{-x},\,0}{pq,\,q^{-N}}\Bigm|q\,;q\Bigr).
\end{align}

\bigskip
\subsection{$q$-Krawtchouk}
The birth and death rates are linear in $q^x$:
\begin{align}
  B(x)&=q^{x-N}-1,\qquad
  D(x)=p(1-q^x),\quad p>0,\\
  \mathcal{E}(n)&=(q^{-n}-1)(1+pq^n),\qquad
  \eta(x)=q^{-x}-1,\\
  \phi_0(x)^2&=\frac{(q\,;q)_N}{(q\,;q)_x(q\,;q)_{N-x}}\,
  p^{-x}q^{\frac12x(x-1)-xN},\\
  d_n^2
  &=\frac{(q\,;q)_N}{(q;q)_n(q;q)_{N-n}}\,
  \frac{(-p\,;q)_n}{(-pq^{N+1}\,;q)_n\,p^nq^{\frac12n(n+1)}}\,
  \frac{1+pq^{2n}}{1+p}
  \times\frac{p^{N}q^{\frac12N(N+1)}}{(-pq\,;q)_N},\\[4pt]
 \check{P}_n(x)&=P_n(\eta(x))
  ={}_3\phi_2\Bigl(
  \genfrac{}{}{0pt}{}{q^{-n},\,q^{-x},\,-pq^n}{q^{-N},\,0}\Bigm|q\,;q\Bigr).
\end{align}

\subsection{quantum $q$-Krawtchouk}
The birth and death rates are quadratic  in $q^x$:
\begin{align}
  B(x)&=p^{-1}q^x(q^{x-N}-1),\qquad
  D(x)=(1-q^x)(1-p^{-1}q^{x-N-1}),\\
  \mathcal{E}(n)&=1-q^n,\qquad
  \eta(x)=q^{-x}-1,\quad p>q^{-N},\\
  \phi_0(x)^2
  &=\frac{(q\,;q)_N}{(q\,;q)_x(q\,;q)_{N-x}}\,
  \frac{p^{-x}q^{x(x-1-N)}}{(p^{-1}q^{-N}\,;q)_x}\,,\\[4pt]
  d_n^2
  &=\frac{(q\,;q)_N}{(q\,;q)_n(q\,;q)_{N-n}}\,
  \frac{p^{-n}q^{-Nn}}{(p^{-1}q^{-n}\,;q)_n}\,
  \times(p^{-1}q^{-N}\,;q)_N,\\[4pt]
 \check{P}_n(x)&= P_n(\eta(x))
  ={}_2\phi_1\Bigl(
  \genfrac{}{}{0pt}{}{q^{-n},\,q^{-x}}{q^{-N}}\Bigm|q\,;pq^{n+1}\Bigr).
\end{align}

\subsection{$q$-Hahn}
The birth and death rates are  quadratic polynomials in $q^x$:
\begin{align}
B(x)&=(1-aq^x)(q^{x-N}-1),\quad
  D(x)= aq^{-1}(1-q^x)(q^{x-N}-b),\quad 0<a,b<1,\\
 \mathcal{E}(n)
  &=(q^{-n}-1)(1-abq^{n-1}),\qquad
  \eta(x)=q^{-x}-1,\\
    \phi_0(x)^2
  &=\frac{(q\,;q)_N}{(q\,;q)_x\,(q\,;q)_{N-x}}\,
  \frac{(a;q)_x\,(b\,;q)_{N-x}}{(b\,;q)_N\,a^x}\,,\\[4pt]
  d_n^2
  &=\frac{(q\,;q)_N}{(q\,;q)_n\,(q\,;q)_{N-n}}\,
  \frac{(a,abq^{-1};q)_n}{(abq^N,b\,;q)_n\,a^n}\,
  \frac{1-abq^{2n-1}}{1-abq^{-1}}
  \times\frac{(b\,;q)_N\,a^N}{(ab\,;q)_N},\\[4pt]
 \check{P}_n(x)&= P_n(\eta(x))
  ={}_3\phi_2\Bigl(
  \genfrac{}{}{0pt}{}{q^{-n},\,abq^{n-1},\,q^{-x}}
  {a,\,q^{-N}}\Bigm|q\,;q\Bigr).
\end{align}

\subsection{dual $q$-Hahn}
The birth and death rates are rational functions of $q^x$:
\begin{align}
  B(x)&=
  \frac{(q^{x-N}-1)(1-aq^x)(1-abq^{x-1})}
  {(1-abq^{2x-1})(1-abq^{2x})},\qquad 0<a,b<1,\\[4pt]
  D(x)&=aq^{x-N-1}
  \frac{(1-q^x)(1-abq^{x+N-1})(1-bq^{x-1})}
  {(1-abq^{2x-2})(1-abq^{2x-1})},\\[4pt]
  \mathcal{E}(n)&=q^{-n}-1,\qquad
  \eta(x)=(q^{-x}-1)(1-abq^{x-1}),\\[4pt]
  \phi_0(x)^2
  &=\frac{(q\,;q)_N}{(q\,;q)_x\,(q\,;q)_{N-x}}\,
  \frac{(a,abq^{-1}\,;q)_x}{(abq^N,b\,;q)_x\,a^x}\,
  \frac{1-abq^{2x-1}}{1-abq^{-1}}\,,\\[4pt]
  d_n^2
  &=\frac{(q\,;q)_N}{(q\,;q)_n\,(q\,;q)_{N-n}}\,
  \frac{(a\,;q)_n(b\,;q)_{N-n}}{(b;q)_N\,a^n}
  \times\frac{(b\,;q)_N\,a^N}{(ab;q)_N}\,,\\[4pt]
\check{P}_n(x)&=P_n(\eta(x))
  ={}_3\phi_2\Bigl(
  \genfrac{}{}{0pt}{}{q^{-n},\,abq^{x-1},\,q^{-x}}
  {a,\,q^{-N}}\Bigm|q\,;q\Bigr).
 \end{align}

\subsection{$q$-Racah}
\label{qrac}
The birth and death rates are rational functions of $q^x$,
\begin{align}
 &B(x)
  =-\frac{(1-aq^x)(1-bq^x)(1-cq^x)(1-dq^x)}
  {(1-dq^{2x})(1-dq^{2x+1})}\,,\\[4pt]
  &D(x)
  =- \tilde{d}\,
  \frac{(1-a^{-1}dq^x)(1-b^{-1}dq^x)(1-c^{-1}dq^x)(1-q^x)}
  {(1-dq^{2x-1})(1-dq^{2x})},
  \label{qracahbd}\\
&c=q^{-N},\ \ a\leq b,\ \ 0<d<1,\ \ 0<a<q^Nd,\  \ qd<b<1, \\
  &\mathcal{E}(n)=(q^{-n}-1)(1-\tilde{d}q^n),\qquad
  \eta(x)=(q^{-x}-1)(1-dq^x),\quad  \tilde{d}\eqdef abcd^{-1}q^{-1},\\
 & \phi_0(x)^2=\frac{(a,b,c,d\,;q)_x}
  {(a^{-1}dq,b^{-1}dq,c^{-1}dq,q\,;q)_x\,\tilde{d}^x}\,
  \frac{1-dq^{2x}}{1-d},\\[4pt]
   &d_n^2
  =\frac{(a,b,c,\tilde{d}\,;q)_n}
  {(a^{-1}\tilde{d}q,b^{-1}\tilde{d}q,c^{-1}\tilde{d}q,q\,;q)_n\,d^n}\,
  \frac{1-\tilde{d}q^{2n}}{1-\tilde{d}}
  \times
  \frac{(-1)^N(a^{-1}dq,b^{-1}dq,c^{-1}dq\,;q)_N\,\tilde{d}^Nq^{\frac12N(N+1)}}
  {(\tilde{d}q\,;q)_N(dq\,;q)_{2N}},\\
&\check{P}_n(x)=
  P_n(\eta(x))
  ={}_4\phi_3\Bigl(
  \genfrac{}{}{0pt}{}{q^{-n},\,\tilde{d}q^n,\,q^{-x},\,dq^x}
  {a,\,b,\,c}\Bigm|q\,;q\Bigr). 
\end{align}

The BD processes corresponding to the 10 polynomials from Krawtchouk (\S\ref{kra}) to 
$q$-Racah (\S\ref{qrac}) are 
on finite lattices, $\mathcal{X}=\{0,1,\ldots,N\}$.

\subsection{Al-Salam-Carlitz II}
The BD processes corresponding to this polynomial are on a semi-infinite lattice $\mathcal{X}=\mathbb{Z}_{\ge0}$.
The birth and death rates are quadratic in $q^x$, and thus bounded,
\begin{align}
B(x)&=aq^{2x+1},\qquad\quad
  D(x)=(1-q^x)(1-aq^x),\quad 0<a<q^{-1},\\
  \mathcal{E}(n)&=1-q^n,\qquad\quad
  \eta(x)=q^{-x}-1,\\
  \phi_0(x)^2&=\frac{a^xq^{x^2}}{(q,aq\,;q)_x}\,,\quad
  d_n^2
  =\frac{(aq)^n}{(q\,;q)_n}\times(aq\,;q)_{\infty}\,,\\[4pt]
  \check{P}_n(x)&=P_n(\eta(x))
  ={}_2\phi_0\Bigl(
  \genfrac{}{}{0pt}{}{q^{-n},\,q^{-x}}{-}\Bigm|q\,;a^{-1}q^n\Bigr).
  \label{alsalamIInorm}
\end{align}

%
%
\section{Explicit Examples II}
\label{ex2}
The solutions of the discrete time BD processes described in this section take different forms from those 
presented in {\bf Theorems \ref{theo3}-\ref{theo4}} \eqref{ctbdsol1}--\eqref{dtbdsol2}, 
although the structure of the problem setting itself \eqref{jacobi}--\eqref{bdcond1} is common.
The eigenvectors of $L$ \eqref{Lndef1} consist of two sets of polynomials $\{\hat{\phi}_n(x)\}$ and 
$\{\hat{\phi}_n^{(-)}(x)\}$, having the same name. They are mutually orthogonal and constitute a
complete set of basis of the corresponding Hilbert space.
These discrete time BD processes are all Markov chains on a semi-infinite lattice 
$\mathcal{X}=\mathbb{Z}_{\ge0}$.
\subsection{$q$-Meixner}
\label{qMei}
The birth and death rates are
quadratic in $q^x$,
\begin{align}
  B(x)&=cq^x(1-bq^{x+1}),\quad
  D(x)=(1-q^x)(1+bcq^x),\quad 0<b<q^{-1},\quad c>0,\\
  \mathcal{E}(n)&=1-q^n,\qquad
  \eta(x)=q^{-x}-1,\\
\phi_0(x)^2&=
  \frac{(bq\,;q)_x}{(q,-bcq\,;q)_x}\,c^xq^{\frac12x(x-1)},\quad
  d_n^2
  =\frac{q^n(bq\,;q)_n}{(q,-c^{-1}q\,;q)_n}
  \times\frac{(-bcq\,;q)_{\infty}}{(-c\,;q)_{\infty}}\,,\\[4pt]
   \check{P}_n(x)&=P_n(\eta(x))
  ={}_2\phi_1\Bigl(
  \genfrac{}{}{0pt}{}{q^{-n},\,q^{-x}}{bq}\Bigm|q\,;-c^{-1}q^{n+1}\Bigr).
\end{align}
It turned out that the above $q$-Meixner polynomials did not form a complete set 
\cite{atakishi1}.
It can be seen clearly by (3.15) of \cite{os34}.
For the completeness another set of orthogonal polynomials is necessary.
They are obtained from the original set by the parameter change (involution)
\begin{align}
  &\hspace{50mm} (b,c)\to(-bc,c^{-1}),\n
  &B^{(-)}(x)= c^{-1}q^x(1+bcq^{x+1}),\quad
  D^{(-)}(x)= (1-q^x)(1-bq^x),
  \label{qM:BD-}\\
 & \check{P}^{(-)}_n(x)
  = P^{(-)}_n\bigl(\eta(x)\bigr)
  ={}_2\phi_1\Bigl(\genfrac{}{}{0pt}{}
  {q^{-n},q^{-x}}{-bcq}\!\!\Bigm|\!q\,;-cq^{n+1}\Bigr),\quad
  \mathcal{E}^{\prime}(n)\eqdef 1+cq^n, \quad \eta(x)=q^{-x}-1,
  \label{qM:Pnm}\\
    &\phi^{(-)}_0(x)^2
  =\frac{(-bcq;q)_x}{(q,bq;q)_x}c^{-x}q^{\frac12x(x-1)},\quad
    \phi^{(-)}_0(x)=(-1)^x\prod_{y=0}^{x-1}
  \sqrt{\frac{B^{(-)}(y)}{D^{(-)}(y+1)}},\quad
  (-1)^x\phi^{(-)}_0(x)>0.
  \label{qMsup}
\end{align}
The orthogonality relations are ($n,m=0,1,\ldots$)
\begin{align}
  (\phi_n,\phi_m)&=\sum_{x=0}^{\infty}
  \phi_n(x)\phi_m(x)
  =\frac{\delta_{n\,m}}{d_n^2},\qquad \quad \hat{\phi}_n(x)\eqdef \phi_n(x)d_n,
  \label{qM:orthorel0}\\
  \bigl(\phi^{(-)}_n,\phi^{(-)}_m\bigr)&=\sum_{x=0}^{\infty}
  \phi^{(-)}_n(x)\phi^{(-)}_m(x)
  =\frac{\delta_{n\,m}}{d^{(-)2}_n},
  \ \ \, d^{(-)}_n\eqdef
  d_n\bigl|_{(b,c)\to(-bc,c^{-1})},\n
  \bigl(\phi_n,\phi^{(-)}_m\bigr)&=\sum_{x=0}^{\infty}
  \phi_n(x)\phi^{(-)}_m(x)=0,\qquad \qquad \hat{\phi}^{(-)}_n(x)\eqdef \phi_n^{(-)}(x)d_n^{(-)}.
  \label{qM:orthorel2}
\end{align}
The  formulas for {\em continuous time} BD involving two sets of orthogonal polynomials are slightly different from the
previous case as given in the following
\begin{theo}
\label{theo51}
The solution  of the initial value problem is
\begin{align} 
\mathcal{P}(x;t)=\hat{\phi}_0(x)\sum_{n=0}^\infty\left(c_ne^{-\mathcal{E}(n)t}\hat{\phi}_n(x)
+c_n^{(-)}e^{-\mathcal{E}^\prime(n)t}\hat{\phi}_n^{(-)}(x)\right),
\label{cBDsolqM}
\end{align}
in which
\begin{align} 
c_n\eqdef\sum_{x=0}^\infty\hat{\phi}_n(x)\hat{\phi}_0(x)^{-1}\mathcal{P}(x;0),
\quad c_n^{(-)}\eqdef\sum_{x=0}^\infty\hat{\phi}_n^{(-)}(x)\hat{\phi}_0(x)^{-1}\mathcal{P}(x;0).
\label{cnqM}
\end{align}
The  transition probability from
site $y$ at time $t=0$ to site $x$ at a later time $t$ in the {\em continuous time} BD is
\begin{equation}
  \mathcal{P}(x,y;t)=\hat{\phi}_0(x)
  \sum_{n=0}^{\infty}\Bigl(
  e^{-\mathcal{E}(n)t}\hat{\phi}_n(x)\hat{\phi}_n(y)
  +e^{-\mathcal{E}^{\prime}(n)t}\hat{\phi}^{(-)}_n(x)\hat{\phi}^{(-)}_n(y)
  \Bigr)\hat{\phi}_0(y)^{-1}\ \ (t>0),
  \label{qMdualtranprob}
\end{equation}
\end{theo}
These formulas are reported in \cite{os34} \S6.A (6.20).

The  formulas for  {\em discrete time} BD involving two sets of orthogonal polynomials are
shown in the following
\begin{theo}
\label{theo52}
The solution  of the initial value problem and the transition matrix after $\ell$ step are
\begin{align} 
\mathcal{P}(x;\ell)&=\hat{\phi}_0(x)\sum_{n=0}^\infty\left(c_n\kappa(n)^\ell\hat{\phi}_n(x)
+c_n^{(-)}{\kappa^{(-)}(n)^\ell}\hat{\phi}_n^{(-)}(x)\right),
\label{DTcBDsolqM}\\
  \mathcal{P}(x,y;\ell)&=\hat{\phi}_0(x)
  \sum_{n=0}^{\infty}\Bigl(
  \kappa(n)^\ell\hat{\phi}_n(x)\hat{\phi}_n(y)
  +{\kappa^{(-)}(n)^\ell}\hat{\phi}^{(-)}_n(x)\hat{\phi}^{(-)}_n(y)
  \Bigr)\hat{\phi}_0(y)^{-1}\ \ \ell=1,2\ldots,
  \label{DTqMdualtranprob}
\end{align}
in which $c_n$ and $c_n^{(-)}$ are given in \eqref{cnqM} and
\begin{equation}
\kappa(n)= 1-t_S\cdot\mathcal{E}(n),\quad 
\kappa^{(-)}(n)= 1-t_S\cdot\mathcal{E}^\prime(n).
\end{equation}
\end{theo}

\begin{theo}
\label{theoSpe2}
The spectral representation of the symmetric matrix $\mathcal{H}$ \eqref{Hdef} for the present case reads
\begin{equation}
\mathcal{H}_{x\,y}=\sum_{n=0}^\infty\mathcal{E}(n)\hat{\phi}_n(x)\hat{\phi}_n(y)
+\sum_{n=0}^\infty\mathcal{E}'(n)\hat{\phi}_n^{(-)}(x)\hat{\phi}_n^{(-)}(y).
\label{Hspe2}
\end{equation}
This in turn supplies the spectral representations of $L_{BD}$ and $L$ through \eqref{LBDHrel} and \eqref{eigrel},
\begin{align} 
 {L_{BD}}_{x\,y}&=-\hat{\phi}_0(x)\left(\sum_{n=0}^\infty\mathcal{E}(n)\hat{\phi}_n(x)\hat{\phi}_n(y)
 +\sum_{n=0}^\infty\mathcal{E}'(n)\hat{\phi}_n^{(-)}(x)\hat{\phi}_n^{(-)}(y)\right)
 \hat{\phi}_0(y)^{-1}\\
& =-\pi(x)\sum_{n=0}^\infty\mathcal{E}(n)({d_n^2}/{d_0^2})\check{P}_n(x)\check{P}_n(y)\n
&\hspace{10mm}-(-1)^{x+y}\sum_{n=0}^\infty\mathcal{E}'(n)
\frac{q^{\tfrac12x(x-1)}}{(q\,;q)_x}\frac{c^{-y}(-bcq\,;q)_y}{(bq\,;q)_y}
\left(d_n^{(-)}\right)^2\check{P}_n^{(-)}(x)\check{P}_n^{(-)}(y),\\ 
 L_{x\,y}&=\hat{\phi}_0(x)\left(\sum_{n=0}^\infty\kappa(n)\hat{\phi}_n(x)\hat{\phi}_n(y)
 +\sum_{n=0}^\infty\kappa^{(-)}(n)\hat{\phi}_n^{(-)}(x)\hat{\phi}_n^{(-)}(y)\right)
 \hat{\phi}_0(y)^{-1}\\
& =\pi(x)\sum_{n=0}^\infty\kappa(n)({d_n^2}/{d_0^2})\check{P}_n(x)\check{P}_n(y)\n
&\hspace{5mm}+(-1)^{x+y}\sum_{n=0}^\infty\kappa^{(-)}(n)
\frac{q^{\tfrac12x(x-1)}}{(q\,;q)_x}\frac{c^{-y}(-bcq\,;q)_y}{(bq\,;q)_y}
\left(d_n^{(-)}\right)^2\check{P}_n^{(-)}(x)\check{P}_n^{(-)}(y), 
\end{align}
in which $\pi(x)$ \eqref{pidef} is the stationary distribution.
\end{theo}

The $q$-Meixner polynomials provide another exactly solvable BD 
with birth rate $B^{(-)}(x)$ and death rate $D^{(-)}(x)$ \eqref{qM:BD-}.
Discrete time BD  is also possible as shown in the following
\begin{theo}
\label{theo53}
The formulas for the initial value problem and the transition matrix for the discrete BD
defined by the second set of polynomials are
\begin{align} 
\mathcal{P}(x;\ell)&=\hat{\phi}_0^{(-)}(x)\sum_{n=0}^\infty\left(\bar{c}_n\kappa(n)^\ell\hat{\phi}_n^{(-)}(x)
+c_n^{(+)}{\kappa^{(+)}(n)}^\ell\hat{\phi}_n(x)\right),
\label{DTcBDsolqM-}\\
\bar{c}_n&\eqdef\sum_{x=0}^\infty\hat{\phi}_n^{(-)}(x)\hat{\phi}_0^{(-)}(x)^{-1}\mathcal{P}(x;0),
\quad c_n^{(+)}\eqdef\sum_{x=0}^\infty\hat{\phi}_n(x)\hat{\phi}_0^{(-)}(x)^{-1}\mathcal{P}(x;0).
\label{cnqM-}\\
  \mathcal{P}(x,y;\ell)&=\hat{\phi}_0^{(-)}(x)
  \sum_{n=0}^{\infty}\Bigl(
  \kappa(n)^\ell\hat{\phi}_n^{(-)}(x)\hat{\phi}_n^{(-)}(y)
  +{\kappa^{(+)}(n)}^\ell\hat{\phi}_n(x)\hat{\phi}_n(y)
  \Bigr)\hat{\phi}_0^{(-)}(y)^{-1}\ \ \ell=1,2\ldots,
  \label{DTqMdualtranprob-}\\
 L_{x\,y}&=\hat{\phi}_0^{(-)}(x)\left(\sum_{n=0}^\infty\kappa(n)\hat{\phi}_n^{(-)}(x)\hat{\phi}_n^{(-)}(y)
 +\sum_{n=0}^\infty\kappa^{(+)}(n)\hat{\phi}_n(x)\hat{\phi}_n(y)\right)
 \hat{\phi}_0^{(-)}(y)^{-1},
\end{align}
in which
\begin{equation}
\kappa(n)= 1-t_S\cdot\mathcal{E}(n),\quad 
\kappa^{(+)}(n)= 1-t_S\cdot\mathcal{E}^{(+)\prime}(n),\quad 
\mathcal{E}^{(+)\prime}(n)\eqdef 1+c^{-1}q^n.
\end{equation}
\end{theo}
The corresponding formulas for the continuous time BD are obtained from those in {\bf Theorem \ref{theo51}}
by changing $\hat{\phi}_n\leftrightarrow \hat{\phi}_n^{(-)}$, $c_n\leftrightarrow c_n^{(-)}$, etc.

\subsection{$q$-Charlier}
\label{qcha}
The $q$-Charlier polynomials are obtained from those of $q$-Meixner by
setting $b=0$ and $c=a>0$.
The birth and death rates and the complete set of orthogonal vectors involving the $q$-Charlier
polynomials are
\begin{align}
&B(x)=aq^x,\quad  D(x)=1-q^x,\qquad B^{(-)}(x)=a^{-1}q^x,\quad  D^{(-)}(x)=1-q^x,\\
  &\check{P}_n(x)=P_n\bigl(\eta(x)\bigr)=
  {}_2\phi_1\Bigl(\genfrac{}{}{0pt}{}{q^{-n},\,q^{-x}}{0}\!\Bigm|\!q\,;
  -a^{-1}q^{n+1}\Bigr),
  \quad \phi_0(x)=\sqrt{\frac{a^xq^{\frac12x(x-1)}}{(q;q)_x}}\,,
  \label{qCsup0}\\
  &\check{P}_n^{(-)}(x)=P_n^{(-)}\bigl(\eta(x)\bigr)=
    {}_2\phi_1\Bigl(\genfrac{}{}{0pt}{}{q^{-n},\,q^{-x}}{0}\!\Bigm|\!q\,;
  -aq^{n+1}\Bigr),
  \ \ \phi^{(-)}_0(x)
  =(-1)^x\sqrt{\frac{a^{-x}q^{\frac12x(x-1)}}{(q;q)_x}}\,,
  \label{qCsup}\\
  & \mathcal{E}(n)=1-q^n,\quad  \mathcal{E}^{\prime}(n)\eqdef 1+aq^n,
  \quad  \mathcal{E}^{(+)\prime}(n)\eqdef 1+a^{-1}q^n,\quad \eta(x)=q^{-x}-1,\\ 
&  d_n^2
  =\frac{q^n}{(q,-a^{-1}q\,;q)_n}
  \times\frac{1}{(-a\,;q)_{\infty}},\qquad d^{(-)}_n\eqdef d_n|_{a\to a^{-1}}. 
 \label{qCsupeig}
\end{align}
The orthogonality relations have the same forms as
\eqref{qM:orthorel0}--\eqref{qM:orthorel2}.
The formulas in {\bf Theorem \ref{theo51}--\ref{theo53}} apply for the $q$-Charlier.

\subsection{dual big $q$-Jacobi}
\label{dbqj}
The big $q$-Jacobi polynomials are the most generic member of the family having orthogonality measures
of Jackson integral type. The big $q$-Laguerre, Al-Salam-Carlitz $\I$, discrete $q$-Hermite $\I$,$\II$ 
and $q$-Laguerre belong to this family.
They all have unbounded $B(x)$ and $D(x)$.
This means that the discrete time BD corresponding to these polynomials cannot be constructed, although
the solutions of their {\em continuous time} BD processes show quite interesting features, 
as reported in \cite{os34} ($\II$\S6A). 

Some of the dual polynomials of this family, however, provide exactly solvable discrete time BD processes
as they have bounded $B(x)$ and $D(x)$.
Reflecting the structure of Jackson integrals, the corresponding dual orthogonal polynomials consist of two
components, similar to the $q$-Meixner case presented in \S\ref{qMei}.
The dual polynomials belonging to this family are extensively reported in \cite{atakishi1}--\cite{atakishi2}.

The data for the dual big $q$-Jacobi polynomials are
\begin{align}
  &0<a<q^{-1},\quad 0<b<q^{-1},\quad c<0, \quad \eta(x)\eqdef(q^{-x}-1)(1-abq^{x+1}),
  \label{bqjpararange}\\
  &B(x)
  =-cq^{x+1}\frac{(1-aq^{x+1})(1-abq^{x+1})(1-abc^{-1}q^{x+1})}
  {(1-abq^{2x+1})(1-abq^{2x+2})},
  \label{dbqJB+}\\
  &D(x)
  =aq\frac{(1-q^x)(1-bq^x)(1-cq^x)}{(1-abq^{2x})(1-abq^{2x+1})}, 
  \label{dbqJD+}\\
  &B^{(-)}(x)
  =aq^{x+1}\frac{(1-bq^{x+1})(1-abq^{x+1})(1-cq^{x+1})}
  {(1-abq^{2x+1})(1-abq^{2x+2})},
\label{dbqJB-}\\
  &D^{(-)}(x)
  =-cq\frac{(1-q^x)(1-aq^x)(1-abc^{-1}q^x)}{(1-abq^{2x})(1-abq^{2x+1})}, 
  \label{dbqJD-}\\
   &\check{P}_n(x)
  = P_n\bigl(\eta(x)\bigr)=
  {}_3\phi_2\Bigl(\genfrac{}{}{0pt}{}
  {q^{-n},\,abq^{x+1},\,q^{-x}}{aq,\,abc^{-1}q}\!\Bigm|\!q\,;
  ac^{-1}q^{n+1}\Bigr),
  \label{Px+}\\
  &\check{P}^{(-)}_n(x)
  = P^{(-)}_n\bigl(\eta(x)\bigr)=
  {}_3\phi_2\Bigl(\genfrac{}{}{0pt}{}
  {q^{-n},\,abq^{x+1},\,q^{-x}}{bq,\,cq}\!\Bigm|\!q\,;
  a^{-1}cq^{n+1}\Bigr),
  \label{Px-}\\
  &\phi_0(x)^2=\frac{q^{\frac12x(x-1)}}{(-ac^{-1})^x}\frac{(abc^{-1}q;q)_x}{(cq;q)_x}
\frac{1-abq^{2x+1}}{1-abq^{x+1}}
  \frac{(aq,abq^2;q)_x}{(q,bq;q)_x},    \\
  &\phi_0^{(-)}(x)^2=\frac{q^{\frac12x(x-1)}}{(-a^{-1}c)^x}\frac{(bq;q)_x}{(aq;q)_x}
\frac{1-abq^{2x+1}}{1-abq^{x+1}}
  \frac{(abq^2,cq;q)_x}{(q,abc^{-1}q;q)_x}, 
  \quad(-1)^x\phi^{(-)}_0(x)>0, 
\label{pmpol}  \\
 &d_n^2 = q^{n}
   \frac{(aq,abc^{-1}q;q)_n}{(q,ac^{-1}q;q)_n}\times
  \frac{(bq,cq;q)_{\infty}}
  {(abq^2,a^{-1}c;q)_{\infty}},\\
  &d_n^{(-)2} =q^n 
  \frac{(bq,cq,;q)_n}
  {(q,a^{-1}cq;q)_n}\times \frac{(aq,abc^{-1}q;q)_{\infty}}
  {(abq^2,ac^{-1};q)_{\infty}}.
\end{align}
The following parameter substitution
(involution)
\begin{equation}
  (a,b,c)\leftrightarrow(c,abc^{-1},a),
  \label{bqjsym}
\end{equation}
gives rise to the interchange of the basic and the $(-)$ objects.
The orthogonality relations for  $\{\phi_n(x)\}$ and $\{\phi_n^{(-)}(x)\}$ 
have the same form as  \eqref{qM:orthorel0}--\eqref{qM:orthorel2}. 

Discrete time BD based on $B(x)$, $D(x)$ \eqref{dbqJB+},\eqref{dbqJD+}
has the same formulas as those in {\bf Theorem \ref{theo52}} given for the 
$q$-Meixner systems \eqref{DTcBDsolqM} and \eqref{DTqMdualtranprob}
with the replacements
\begin{equation}
\mathcal{E}(n)= aq(1-q^n), \quad \mathcal{E}^\prime(n)=q(a-cq^n).
\end{equation}
Likewise discrete time BD based on $B^{(-)}(x)$, $D^{(-)}(x)$ \eqref{dbqJB-},\eqref{dbqJD-}
has the same forms as those  in {\bf Theorem \ref{theo53}} \eqref{DTcBDsolqM-}--\eqref{DTqMdualtranprob-}
with the replacements
\begin{equation}
\mathcal{E}(n)= -cq(1-q^n), \quad \mathcal{E}^{(+)\prime}(n)=q(-c+aq^n).
\end{equation}

\subsection{dual big $q$-Laguerre}
\label{dbql}
The basic data are as follows.
They are  obtained from those of the dual big $q$-Jacobi 
polynomial by setting
$b\to0$ and $c\to b$, with
$0<a<q^{-1}$ and $b<0$,
\begin{align}
  B(x)&=-bq^{x+1}(1-aq^{x+1}),
  \ \ D(x)=aq(1-q^x)(1-bq^x),
  \label{dbqlB+}\\
  B^{(-)}(x)&=aq^{x+1}(1-bq^{x+1}),
  \quad \ D^{(-)}(x)=-bq(1-q^x)(1-aq^x),
  \label{dbqlB-}\\
   \check{P}_n(x)&=
  P_n\bigl(\eta(x)\bigr)=
  {}_2\phi_1\Bigl(\genfrac{}{}{0pt}{}{q^{-n},\,q^{-x}}{aq}\!\Bigm|\!q\,;
  ab^{-1}q^{n+1}\Bigr),\qquad \eta(x)\eqdef q^{-x}-1,
  \label{dbqlQx+}\\
  \check{P}^{(-)}_n(x)&=
  P^{(-)}_n\bigl(\eta(x)\bigr)=
  {}_2\phi_1\Bigl(\genfrac{}{}{0pt}{}{q^{-n},\,q^{-x}}{bq}\!\Bigm|\!q\,;
  a^{-1}bq^{n+1}\Bigr),
  \label{bqlQx-}\\
\phi_0(x)^2
 &=\frac{q^{\frac12x(x-1)}}{(-ab^{-1})^x}\frac{1}{(bq;q)_x}
  \frac{(aq;q)_x}{(q;q)_x},  
  \\
\phi_0^{(-)}(x)^2
 &=\frac{q^{\frac12x(x-1)}}{(-a^{-1}b)^x}\frac{1}{(aq;q)_x}
  \frac{(bq;q)_x}{(q;q)_x}, 
  \qquad(-1)^x\phi^{(-)}_0(x)>0,   \\
 d_n^2&= q^{n}
  \frac{(aq;q)_n}{(q,ab^{-1}q;q)_n}
  \times \frac{(bq;q)_{\infty}}{(a^{-1}b;q)_{\infty}},\\
  d_n^{(-)2}&=q^n 
  \frac{(bq;q)_n}{(q,a^{-1}bq;q)_n}
  \times  \frac{(aq;q)_{\infty}}{(ab^{-1};q)_{\infty}},\\
\phi_n^{(-)}(x)&\eqdef \phi_0^{(-)}(x)\check{P}_n^{(-)}(x), \quad \hat{\phi}_n^{(-)}(x)\eqdef \phi_n^{(-)}(x)d_n^{(-)}.
\label{pmpol2}
\end{align}
The basic objects and the $(-)$ objects are interchanged by the parameter substitution (involution) 
$a\leftrightarrow b$.

\bigskip
Discrete time BD based on $B(x)$, $D(x)$ \eqref{dbqlB+}
has the same formulas as those in {\bf Theorem \ref{theo52}} for the $q$-Meixner systems \eqref{DTcBDsolqM} and \eqref{DTqMdualtranprob}
with the replacements
\begin{equation}
\mathcal{E}(n)= aq(1-q^n), \quad \mathcal{E}^\prime(n)=q(a-bq^n).
\end{equation}
Likewise discrete time BD based on $B^{(-)}(x)$, $D^{(-)}(x)$ \eqref{dbqlB-}
has the same forms as those in {\bf Theorem \ref{theo53}}  \eqref{DTcBDsolqM-}--\eqref{DTqMdualtranprob-}
with the replacements
\begin{equation}
\mathcal{E}(n)= -bq(1-q^n), \quad \mathcal{E}^{(+)\prime}(n)=q(-b+aq^n).
\end{equation}

%
%
%
\section{Mirror symmetric Birth and Death processes}
\label{msbd}
Here we present a simple mirror symmetric discrete time BD process.
Among the exactly solvable discrete time BD's on a finite integer lattice  listed in section \ref{ex1}, 
those related to two polynomials, the Krawtchouk \S\ref{kra} and Hahn \S\ref{hah}, can be made
mirror symmetric 
\begin{align}
\text{Mirror symmetry:}\quad &D(N-x)=B(x),\ \quad B(N-x)=D(x),
\label{lrsym1}\\
\Longrightarrow\quad &\phi_0(N-x)^2=\phi_0(x)^2,\quad
\check{P}_n(N-x)=(-1)^n\check{P}_n(x),
\label{lrsym2}
\end{align}
by adjusting the parameters. These two polynomials have $\eta(x)=x$. 
As for the Krawtchouk \eqref{kra1}--\eqref{kra4} with $p=1/2$, we have 
\begin{align}
\text{Krawtchouk:}\quad &B(x)=(N-x)/2, \quad D(x)=x/2,\quad \eta(x)=x,\quad \mathcal{E}(n)= n,\\
  \phi_0(x)^2&=
  \frac{N!}{x!\,(N-x)!},\quad
  d_n^2
  =\frac{N!}{n!\,(N-n)!}\times2^{-N},\\
&\check{P}_n(x)=
    P_n(x)
  ={}_2F_1\Bigl(
  \genfrac{}{}{0pt}{}{-n,\,-x}{-N}\Bigm|2\Bigr)=(-1)^nP_n(N-x),
\label{mirkra}
\end{align}
due to 
Pfaff's transformation formula(see \cite{askey}p79)
\begin{equation*}
{}_2F_1\Bigl(
  \genfrac{}{}{0pt}{}{a,\,b}{c}\Bigm| x\Bigr)=(1-x)^{-a}{}_2F_1\Bigl(
  \genfrac{}{}{0pt}{}{a,\,c-b}{c}\Bigm|\frac{x}{x-1}\Bigr).
\end{equation*}
Taking $a=b$ for the Hahn \eqref{hah1}--\eqref{hah4}, we obtain
\begin{align}
\text{Hahn:}\quad &B(x)=(x+a)(N-x), \quad D(x)=x(a+N-x),\\
 &\eta(x)=x,\quad
\mathcal{E}(n)= n(n+2a-1),\\
\phi_0(x)^2
 &=\frac{N!}{x!\,(N-x)!}\,\frac{(a)_x\,(a)_{N-x}}{(a)_N},\quad
  d_n^2
  =\frac{N!}{n!\,(N-n)!}\,
  \frac{(2n+2a-1)(2a)_N}{(n+2a-1)_{N+1}}
  \times\frac{(a)_N}{(2a)_N},\\[4pt]
 &\check{P}_n(x)= P_n(x)
  ={}_3F_2\Bigl(
  \genfrac{}{}{0pt}{}{-n,\,n+2a-1,\,-x}
  {a,\,-N}\Bigm|1\Bigr)=(-1)^nP_n(N-x).
  \label{mirhah}
\end{align}
due to the following transformation formula (see \cite{askey}p142)
\begin{equation*}
{}_3F_2\Bigl(
  \genfrac{}{}{0pt}{}{-n,\,a,\,b}
  {d,\,e}\Bigm|1\Bigr)=\frac{(e-a)_n}{(e)_n}{}_3F_2\Bigl(
  \genfrac{}{}{0pt}{}{-n,\,a,\,d-b}
  {d,\,a+1-n-e}\Bigm|1\Bigr).
\end{equation*}

Let us consider the following transition matrix $L_{BD}^M$
\begin{align}
(L_{BD}^M\mathcal{P})(x;t)&=\sum_{y\in\mathcal{X}}{L_{BD}^M}_{x\,y}\mathcal{P}(y;t)
\label{bdeqformalM}\\
&=-(B(x)+D(x))\mathcal{P}(N-x;t)+B(x-1)\mathcal{P}(N-x+1;t)\n
&\hspace{56.5mm}+D(x+1)\mathcal{P}(N-x-1;t),
\label{BDeqM}
\end{align}
\begin{align*}
&L_{BD}^M\\
&=\left(
\begin{array}{cccccc}
0 &    0& \cdots&\cdots&D(1)&-B(0) \\
0 &  0& \cdots&D(2) &-B(1)-D(1)&B(0) \\
0  &\cdots &\cdots  &  -B(2)-D(2)&B(1) &0\\
\vdots&\cdots&\cdots&\cdots&\cdots&\vdots\\
0&D(N\!-\!1)&\cdots&\cdots&\cdots&0\\
D(N)&-B(N\!-\!1)-D(N\!-\!1)&\cdots&\cdots&0&0\\
-D(N)&B(N\!-\!1)&\cdots&\cdots&0&0
\end{array}
\right),
\end{align*}
which is obtained from $L_{BD}$ \eqref{LBDdef} by mirror reflection.
That is 
\begin{align} 
&\hspace{30mm} L_{BD}^M=L_{BD}J,\\
 J&\eqdef \text{anti-diagonal}\{1,1,\ldots,1\},\ \text{or}\ J_{x\,y}\eqdef \delta_{x,\,N-y},
\quad \sum_{y\in\mathcal{X}}J_{x\,y}\check{P}_n(y)=(-1)^n\check{P}_n(x).
\label{Jdef}
\end{align}
From the general relation \eqref{LBDsol}
\begin{equation*}
(L_{BD}\hat{\phi}_0\hat{\phi}_n)(x)=-\mathcal{E}(n)\hat{\phi}_0(x)\hat{\phi}_n(x),\quad
n\in\mathcal{X},
\end{equation*}
we obtain
\begin{equation}
(L_{BD}^M\hat{\phi}_0\hat{\phi}_n)(x)=-(-1)^n\mathcal{E}(n)\hat{\phi}_0(x)\hat{\phi}_n(x),\quad
n\in\mathcal{X}.
\end{equation}
Since the spectrum of $L_{BD}^M$ is not negative semi-definite, $L_{BD}^M$ does not define
a stochastic process.
However, it is interesting to consider its discrete time version.

\bigskip
Let us introduce the mirror image of $L$ \eqref{Lndef1},
\begin{align}
&L^{M}_{x+1\,N-x}=\bar{B}(x),\ \  L^{M}_{x-1\,N-x}=\bar{D}(x),\quad
 \ L^{M}_{x\,N-x}=1-\bar{B}(x)-\bar{D}(x), \n 
& \qquad L^{M}_{x\,y}=0 \quad \text{for} \quad x+y<N-1,\quad x+y>N+1,
\label{MLndef1}
\end{align}
\begin{align*}
&L^{M}\eqdef LJ=\n
&{\footnotesize
\left(
\begin{array}{cccccc}
0 &    0& \cdots&\cdots&\bar{D}(1)&1-\bar{B}(0) \\
0 &  0& \cdots&\bar{D}(2) &1-\bar{B}(1)-\bar{D}(1)&\bar{B}(0) \\
0  &\cdots &\cdots  &  1-\bar{B}(2)-\bar{D}(2)&\bar{B}(1) &0\\
\vdots&\cdots&\cdots&\cdots&\cdots&\vdots\\
0&\bar{D}(N\!-\!1)&\cdots&\cdots&\cdots&0\\
\bar{D}(N)&1-\bar{B}(N\!-\!1)-\bar{D}(N\!-\!1)&\cdots&\cdots&0&0\\
1-\bar{D}(N)&\bar{B}(N\!-\!1)&\cdots&\cdots&0&0
\end{array}
\right).
}
\end{align*}
This defines another exactly solvable stochastic process, as $L^M$ is a non-negative tri-anti-diagonal
matrix with
\begin{equation}
\left(L^{M}\hat{\phi}_0\hat{\phi}_n\right)(x)=\kappa_{M}(n)\hat{\phi}_0(x)\hat{\phi}_n(x),\quad 1\ge \kappa_{M}(n)\eqdef (-1)^n\left(1-t_S\cdot\mathcal{E}(n)\right)\ge-1,\quad n\in\mathcal{X}.
\label{eigvalM}
\end{equation}
However, this is simply the mirror image of the original process.
Let us introduce the process governed by the sum of $L$ and $L^M$,
\begin{equation}
L^{S}\eqdef\frac12(L+L^{M})=\frac12L(I_d+J),
\label{Slndef1}
\end{equation}
which is exactly solvable having interesting properties,
\begin{align}
&\left(L^{S}\hat{\phi}_0\hat{\phi}_n\right)(x)=\kappa_{S}(n)\hat{\phi}_0(x)\hat{\phi}_n(x),\\
&1\ge \kappa_{S}(n)\eqdef \frac12\left(1+(-1)^n\right)\bigl(1-t_S\cdot\mathcal{E}(n)\bigr)\ge-1,\quad n\in\mathcal{X}.
\label{eigvalS}
\end{align}
All odd eigenvalues vanish 
\begin{equation*}
\kappa_{S}(2n+1)=0,\quad n=0,1,\ldots,[ \tfrac{N-1}2],
\end{equation*}
and the convergence to the stationary distribution is accelerated.
By expanding the initial distribution as \eqref{cndef}, the distribution after $\ell$ step is
described by the even eigenvectors only
\begin{align}
\mathcal{P}(x;\ell)=\hat{\phi}_0(x)\sum_{n=0}^{[\frac{N}2]}c_{2n}\bigl(\kappa_{S}(2n)\bigr)^{\ell}\hat{\phi}_{2n}(x),\quad \ell=1,2,\ldots.
\end{align}
Likewise, the  transition matrix from $y$ at $\ell=0$ ($\mathcal{P}(x;0)=\delta_{x,y}$) to $x$ after $\ell$ steps is
\begin{equation}
\mathcal{P}(x,y;\ell)=\hat{\phi}_0(x)\sum_{n=0}^{[\frac{N}2]}\bigl(\kappa_{S}(2n)\bigr)^\ell\hat{\phi}_{2n}(x)\hat{\phi}_{2n}(y)\hat{\phi}_0(y)^{-1}.
\label{dtbdsol2S}
\end{equation}
In other words, the asymmetric part of the initial distribution $\mathcal{P}^{AS}(x;0)$
\begin{align*}
&\mathcal{P}^{AS}(x;0)\eqdef\frac12\left(\mathcal{P}(x;0)-\mathcal{P}(N-x;0)\right)=\left(\frac12(I_d-J)\mathcal{P}\right)(x;0),
\end{align*}
is erased by one action of $L^S$ \eqref{Slndef1}
\begin{equation*}
\left(L^{S}\mathcal{P}^{AS}\right)(x;0)=\frac14\left(L(I_d+J)(I_d-J)\right)\mathcal{P}(x;0)=0.
\end{equation*}
A Markov chain having similar eigenvectors was reported in \cite{diaconis20}.
A very special case of discrete time BD with $\bar{B}(x)+\bar{D}(x)=1$ based on $p=1/2$ 
Krawtchouk was reported in \cite{hokonno}.

\subsection{Dual systems}
For finite systems, an apparently different looking exactly solvable system can be constructed from a
known exactly solvable one by the similarity transformation by the anti-diagonal matrix $J$ \eqref{Jdef},
\begin{equation*}
L^d\eqdef JLJ,\qquad L^d_{BD}\eqdef JL_{BD}J.
\end{equation*}
That is,
\begin{equation*}
\bar{B}^d(x)\eqdef\bar{D}(N-x),\quad \bar{D}^d(x)\eqdef\bar{B}(N-x);\qquad 
B^d(x)\eqdef{D}(N-x),\quad {D}^d(x)\eqdef {B}(N-x),
\end{equation*}
and it is  called a `dual system' \cite{coo-hoa-rah77}\S7.
The original and its dual system have common eigenvalues and the eigenvectors are mapped by $J$.
For most of examples listed in \S\ref{ex1} dual polynomials are the same as the original one with 
parameter cahnge (involution). 
\section{Comments}
It should be stressed that the same solution procedures apply to BD processes related to 
various new orthogonal polynomials \cite{os26}. 
They are obtained from the classical orthogonal polynomials, {\em e.g.} the Racah and $q$-Racah, etc
 by multiple applications of the discrete analogue of the Darboux transformations or 
the Krein-Adler transformations. Since there are many different ways to deform 
the classical orthogonal polynomials, these new polynomials offer virtually infinite examples of 
exactly solvable BD processes.

%

\end{document}